\documentclass[aop,seceqn,citesort,MSNbibl,dvips]{arximspdf}
\usepackage{graphics}
\doi{10.1214/09-AOP514}
\volume{38}
\issue{4}
\pubyear{2010}
\firstpage{1507}
\lastpage{1531}

\makeatletter

\newtheorem{theorem}{Theorem}[section]
\newtheorem{lemma}[theorem]{Lemma}
\newtheorem{claim}[theorem]{Claim}

\newproclaim{Remark}{Remark}

\makeatother

\begin{document}
\begin{frontmatter}

\title{Poisson polytopes}
\runtitle{Poisson polytopes}

\begin{aug}
\author[A]{\fnms{Imre} \snm{B\'ar\'any}\thanksref{t1,t2}\ead[label=e1]{barany@renyi.hu}} and
\author[B]{\fnms{Matthias} \snm{Reitzner}\corref{}\thanksref{t1}\ead[label=e2]{matthias.reitzner@uni-osnabrueck.de}}
\runauthor{I. B\'ar\'any and M. Reitzner}
\affiliation{R\'enyi Institute of Mathematics and
University of Osnabr\"uck}
\address[A]{R\'enyi Institute of Mathematics\\
Hungarian Academy of Sciences\\
PO Box 127\\
1364 Budapest\\
Hungary\\
and\\
Department of Mathematics\\
University College London\\
Gower Street\\
London WC1E 6BT\\
United Kingdom\\
\printead{e1}} 
\address[B]{Institute of Mathematics\\
University of Osnabr\"uck\\
49069 Osnabr\"uck\\
Germany\\
\printead{e2}}
\end{aug}

\thankstext{t1}{Supported in part by the European
Network PHD, MCRN-511953.}

\thankstext{t2}{Supported by Hungarian National
Foundation Grants T 60427 and T 62321.}

\pdfauthor{Imre Barany, Matthias Reitzner}

\received{\smonth{9} \syear{2008}}
\revised{\smonth{9} \syear{2009}}

%
\begin{abstract}
We prove the central limit theorem for the volume and the $f$-vector of
the Poisson random polytope $\Pi_{\eta}$ in a fixed convex polytope $P
\subset\mathbb{R}^d$. Here, $\Pi_{\eta}$ is the convex hull of the
intersection of a Poisson process $X$ of intensity $\eta$ with $P$.
\end{abstract}

%

\end{frontmatter}

\section{Introduction and main results}\label{sec:intro}

Let $K\subset\mathbb{R}^d$ be a convex set of volume $1$. Assume that
$X=X(\eta)$ is a Poisson point process in $\mathbb{R}^d$ of intensity
$\eta$.
The intersection of $K$ with $X(\eta)$ consists of uniformly
distributed random points $X_1, \ldots,X_N$ (where $N$ is a random variable).
Define the \textit{Poisson polytope} $\Pi_\eta$, as the convex hull
$ [X_1, \ldots,X_N] = [K \cap X(\eta) ]$.

The study of properties of random convex hulls is a classical subject
in stochastic geometry and
dates back to 1864. Due to the geometric nature of the available methods,
for over one hundred years, investigations
mainly concentrated on the expectation of functionals of random convex hulls
such as volume or number of vertices;
see, for example, the survey of Weil and Wieacker \cite{ww}.

The first distributional results were only proven twenty years ago. In
1988, Groeneboom \cite
{Gro} obtained the central limit theorem (CLT) for the number of
vertices of the
Poisson polytope when the convex body $K$ is the planar disc. In 1994,
a~CLT for the area
of a random polygon in the planar disc was proven by Hsing~\cite{Hsi}.
Recently, this was
generalized to arbitrary dimensions by Reitzner \cite{Re8}, who
established a CLT for
$V(\Pi_{\eta})$, the volume of the Poisson polytope, and for $f_{\ell
}(\Pi_{\eta})$, the
number of $\ell$-dimensional faces of the Poisson polytope, when the
body $K \subset\mathbb{R}^d$
has smooth boundary.

The situation seems to be much more involved when the underlying convex
set is a polytope
$P$. In the planar case, when $P$ is a convex polygon, a CLT for the
number of vertices
$f_0(\Pi_\eta)$ was proven by Groeneboom \cite{Gro} and a CLT for the
area of $\Pi_{\eta}$
by Cabo and Groeneboom \cite{CaGro}, but it seems that the stated
variances are incorrect
(see the discussion in Buchta \cite{Bu11}).

The main result of the present paper is the central limit theorem for
the Poisson polytope
$\Pi_\eta$ for all dimensions $d \geq2$, when the mother body is a
polytope in $\mathbb{R}^d$.
\begin{theorem}\label{th:CLTPi} There exists a function ${\varepsilon
}(\eta)$,
tending to zero as $\eta\to\infty$, such that
for every polytope $P\subset\mathbb{R}^d$ of volume $1$,
\[
\sup_{x \in\mathbb{R}} \biggl\vert\mathbb{P}\biggl( \frac{V(\Pi_{\eta}) -
\mathbb{E} V(\Pi_{\eta})}
{\sqrt{\operatorname{Var}
V(\Pi_{\eta})} }\leq x \biggr) - \Phi(x) \biggr\vert\leq c (P)
{\varepsilon}(\eta)
\]
and for all $\ell=0,\ldots,d-1$,
\[
\sup_{x \in\mathbb{R}} \biggl\vert\mathbb{P}\biggl( \frac{f_{\ell}(\Pi
_{\eta}) - \mathbb{E}
f_{\ell
}(\Pi_{\eta})}
{\sqrt{\operatorname{Var}f_{\ell}(\Pi_{\eta})} }\leq x \biggr) - \Phi
(x) \biggr\vert
\leq c (P) {\varepsilon}(\eta),
\]
where $c(P)$ is a constant depending only on $P$.
\end{theorem}
\begin{Remark*} It will turn out that the error term in Theorem
\ref{th:CLTPi} is
\[
{\varepsilon}(\eta)= (\ln\eta)^{-({d-1})/2 +o(1)}.
\]
The constant $c(P)$ depends on the dimension and a power of
$F(P)$, the number of flags of the polytope $P$. A \textit{flag} is a sequence
of faces $F_0,F_1,\ldots,F_{d-1}$ of $P$ such that for all $i$,
$\dim F_i=i$ and $F_i \subset F_{i+1}$.
\end{Remark*}

The Poisson polytope $\Pi_\eta$ is closely related to the random
polytope $P_{n}$ defined in
the following way: fix $n \in\mathbb{N}$ and choose $n$ random points
$X_1,\ldots,X_{n}$
independently and uniformly from $K$. The \textit{random polytope} $P_{n}$
is just the convex
hull of these points: $P_{n}=[X_1,\ldots,X_{n}]$. Clearly, $P_n$ equals
in distribution the
Poisson polytope $\Pi_\eta$, given that the (Poisson-distributed)
number of points of $X \cap
K$ is precisely $n$.

Starting with R\'enyi and Sulanke \cite{rs1} in 1963, there have been
many results
concerning various properties of $P_{n}$ as $n \to\infty$. For
instance, the
asymptotic behavior of the expectation of the volume $V(P_{n})$, and of
the number,
$f_{\ell}(P_{n})$, of $\ell$-dimensional faces of $P_{n}$
($\ell=0,\ldots,d-1$), have been determined as $n \to\infty$; see
\cite
{ww} for an
extensive survey and also \cite{Bar7,Bar8} and \cite{Re7} for more
recent results.
These results on $P_{n}$ imply immediately analogous results for the
Poisson polytope~$\Pi_\eta$.
For the sake of completeness, we state here the results concerning the
expected volume and number of faces.
\begin{theorem}\label{th:EVPi} Assume that $P$ is a polytope of volume
$1$. Then
\begin{eqnarray*}
1-\mathbb{E} V(\Pi_\eta) &=& \frac{F(P)}{(d+1)^{d-1}(d-1)!} \eta
^{-1} \ln^{d-1}
\eta\bigl(1+o(1)\bigr),
\\
\mathbb{E} f_{\ell}(\Pi_\eta) &=& c(d,\ell)F(P) \ln^{d-1} \eta\bigl(1+o(1)\bigr),
\end{eqnarray*}
where $c(d,\ell)>0$ is a constant depending on $d$ and $\ell$.
\end{theorem}

Somewhat surprisingly, the values of these
expectations are not needed for the proofs of our main theorems.

The proof of Theorem \ref{th:CLTPi} is not simple. It uses a
combination of ideas from probability theory and convex geometry.
Section \ref{sec:plan} contains a short sketch of this proof. First,
we have to introduce some notation and background, and, more
importantly, the economic cap covering theorem that will be used
repeatedly. This is the content of the next section.

\section{Notation and background}\label{sec:notation}

The unit sphere is denoted $S^{d-1}$. As usual, $h_K(u)$ denotes the
support function of $K$ in direction $u\in S^{d-1}$:
\[
h_K(u)=\max\{u\cdot x\dvtx x \in K\}.
\]

A \textit{cap} $C$ of $K$ is the intersection of $K$ with a closed half-space.
This half-space can be written as $ \{ x\in\mathbb{R}^d \vert u \cdot x
\geq h_K(u) -t \}$ with $ u \in S^{d-1}$. Thus,
\[
C = K \cap\{ x\in\mathbb{R}^d \vert u \cdot x \geq h_K(u) -t \}.
\]
The \textit{bounding hyperplane} of $C$ is the one with equation $u
\cdot x = h_K(u) -t$. We define, for ${\lambda}>0$, $C^{\lambda}$ by
\[
C^{\lambda}=K \cap\{ x\in\mathbb{R}^d \vert u \cdot x \geq h_K(u) -
{\lambda}t \}.
\]

An important role throughout is played by the function $v \dvtx K \to
\mathbb{R}
$, defined as
\[
v(z)=\min\{V(K\cap H)\dvtx H \mbox{ is a half-space and } z \in H\}.
\]
The \textit{floating body} with parameter $t$ is just the level set
$K(v\ge t)=\{z \in K\dvtx v(z)
\geq t\}$, which is clearly convex. The \textit{wet part} is $K(v\le t)$,
that is, where $v$ is
at most~$t$. The name comes from the three-dimensional picture where
$K$ is
a container containing
$t$ units of water.

The \textit{minimal cap} of $z \in K$ is a cap $C(z)= C_K (z)$
containing $z$ such that $v(z)=V(C(z))$. It need not be unique. The
\textit{center} of the cap $C = K \cap\{ x\in\mathbb{R}^d \dvtx u
\cdot x
\geq
h_K(u) -t \}$ is a point $x \in\partial K$ with $u \cdot x =
h_K(u)$. The center, again, need not be unique, but this will cause
no harm. Assuming that $x$ is the center of $C $, observe that, for
${\lambda}\geq1$,
\[
C^{\lambda}\subset x + {\lambda}(C-x)
\]
and thus $V(C^{\lambda}) \leq{\lambda}^ d V(C)$ always holds. Also,
$\frac{{\lambda}} d
V(C) \leq V(C^{\lambda})$ holds as long as ${\lambda}t$ is smaller
than the width
of $K$ in the direction $u$. The proof is simple: let $L$ be the section
that has maximal $(d-1)$-dimensional volume among all sections of
the form
\[
K \cap\{ x\in\mathbb{R}^d \dvtx u \cdot x \geq h_K(u) -\tau\}
\qquad\mbox{when }
\tau\in[0,t].
\]
Then $V(C) \le t V_{d-1}(L)$. Here, $V_{d-1}$ stands for
$(d-1)$-dimensional volume. On the other hand, the double cone with
base $L$, apexes $x$ and a point in $K \cap\{ x\in\mathbb{R}^d \dvtx u
\cdot x = h_K(u) -{\lambda}t \}$ is contained in $C^{{\lambda}}$ and
its volume
is at least $\frac{{\lambda}t}d V_{d-1}(L)$. So $\frac{{\lambda}} d
V(C) \le
\frac{{\lambda}t}d V_{d-1}(L) \le V(C^{\lambda})$, which is the
inequality we
wanted to prove.

Analogously for $0 < \mu< 1$, we have $\mu^ d V(C) \leq V(C^\mu)
\leq d \mu V(C)$. For the proof, define $D=C^\mu$ and ${\lambda
}=1/\mu>1$.
Then $D$ is a cap of $K$ and $D^{\lambda}=C$. The inequalities $\frac
1 d
{\lambda}V(D) \leq V(D^{\lambda}) \leq{\lambda}^ d V(D)$ translate
directly to $\mu^ d
V(C) \leq V(C^\mu) \leq d \mu V(C)$. These inequalities will be used
often. We call them the \textit{trivial volume estimates}:
\begin{eqnarray*}
\frac{{\lambda}} d V(C) &\le& V(C^{\lambda}) \leq{\lambda}^d V(C)
\qquad\mbox{for } {\lambda}\geq1, \\
\mu^ d V(C) &\leq& V(C^\mu) \leq d \mu V(C)
\qquad\mbox{for } 0 \leq\mu\leq1,
\end{eqnarray*}
where the left-hand side of the first inequality only holds for $C^
{\lambda}\neq K$ and the right-hand side of the second inequality only for
$C \neq K$. The \textit{Macbeath region}, or \textit{$M$-region} for short,
with center $z$ and factor ${\lambda}>0$ is
\[
M(z,{\lambda})=M_K(z,{\lambda})=z+ {\lambda}[(K-z)\cap(z-K)] .
\]
The $M$-region with ${\lambda}=1$ is just the intersection of $K$ and $K$
reflected with respect to $z$. Thus, $M(z,1)$ is convex and centrally
symmetric with center $z$ and $M(z,{\lambda})$ is a homothetic copy of
$M(z,1)$ with center $z$ and factor of homothety ${\lambda}$. We
define the
function $u\dvtx K \to\mathbb{R}$ by
\[
u(z)= V(M(z,1)).
\]

These definitions are from \cite{ELR,BL} and \cite{Bar1}.
The following results come from the same sources. We will use them
extensively. We assume that $K \subset\mathbb{R}^d$ is a convex body
of volume
$1$. Set
%
%
\begin{equation}
\label{def:s0}
s_0=(2d)^{-2d} .
\end{equation}
\begin{lemma}\label{l:A} If $M(x, \frac12)\cap M(y,\frac12)\ne
\varnothing$, then $M(x,1) \subset
M(y,5)$.
\end{lemma}
\begin{lemma}\label{l:B}If $C$ is a cap, $z \in C$ and ${\lambda}>0$, then
$K\cap M(z,{\lambda})
\subset C^{{\lambda}+1}$.
\end{lemma}
\begin{lemma}\label{l:C} If the cap $C$ is contained in the $M$-region
$M(z,\mu)$ and ${\lambda}
>0$, then $C^{\lambda}\subset M(z,{\lambda}\mu)$.
\end{lemma}
\begin{lemma}\label{l:E} If the bounding hyperplane of a cap $C$ is
tangent to $K(v\ge s)$,
then $s \le V(C) \le ds$.
\end{lemma}

Let $K(v=s)= \partial K(v \geq s)$. Assume that $s \le s_0$ and choose a
maximal system of points $Z=\{z_1, \ldots, z_m\}$ on $K(v=s)$ having
pairwise disjoint Macbeath regions $M(z_i, \frac12)$. Such a system
will be called \textit{saturated}. Note that $Z$ (and even $m$) is not
uniquely defined. However, for each $K$ and $s$, we fix a saturated
system $Z$. We write $Z(s)$ and $m(s)=|Z(s)|$ when we want to
emphasize that our fixed saturated system comes from the level set
$K(v=s)$. Clearly, $V(C(z_i))=s$. Set
\[
K'_i (s) = M \bigl(z_i, \tfrac12 \bigr) \cap C(z_i) \quad\mbox{and}\quad K_i (s) = C^{6}
(z_i).
\]
Note that $K_i(s)$ is a cap of $K$ and so, for ${\lambda}>0$, the set
$K_i^{\lambda}(s)=C^{6{\lambda}}(z_i)$ is another cap of $K$.

The sets $K'_i(s)$ and $K_i(s)$ for $i=1,\ldots,m(s)$ form what is
called an \textit{economic cap covering} in the paper of B\'ar\'any and
Larman \cite{BL}. The following result, the \textit{economic cap
covering theorem}, comes from Theorem 6 in \cite{BL} and Theorem 7
in~\cite{Bar1}.
\begin{theorem}\label{th:capcov} For all $s \in(0,s_0]$ and all convex
bodies $K \subset\mathbb{R}^d$ with \mbox{$V(K) =1$}, we have:
\begin{longlist}
\item $ \bigcup_1^{m(s)} K'_i (s) \subset K(v \leq s) \subset
\bigcup_1^{m(s)} K_i(s)$;

\item $s \le V(K_i (s)) \leq6 ^d s$, $i=1, \ldots, m(s)$;

\item $(6d)^{-d}s \le V(K'_i(s)) \leq2^{-d}s$, $i=1, \ldots, m(s)$;

\item every $C$ with $V(C)\le s$ is contained in $M(z_i,15d)
\subset K^{3d}_i(s)$ for some $i$.
\end{longlist}
\end{theorem}

The sets $K_i'(s)$ are pairwise disjoint, all of them have volume
$\ge(6d)^{-d}s$ and are all contained in $K(v\le s)$. This gives
an upper bound for $m(s)$. Similarly, the sets $K_i(s)$ cover
$K(v\le s)$ and all of them have volume $\le6^ds$. This gives a lower
bound for $m(s)$. These simple arguments will be used repeatedly
and we call them the \textit{usual volume arguments}. Summarizing, we
have
%
%
\begin{equation} \label{ms}
\frac1 {6^d s} V\bigl(K(v \leq s)\bigr) \leq m(s)
\leq\frac1 {(6d)^{-d} s} V\bigl(K(v \leq s)\bigr)
\end{equation}
for $s \leq s_0$.

The economic cap covering theorem has the following direct
consequence.
\begin{claim}\label{cl:lambda}For $s \leq s_0$ and $ {\lambda}>1$,
\[
K( v \leq{\lambda}s) \subset\bigcup K^{3 d^2 {\lambda}}_i (s) .
\]
\end{claim}
\begin{pf}
It is clear that $K( v \leq{\lambda}s)$ is contained
in the
union of all caps $C$
with $V(C)={\lambda}s$. Let $C$ be a cap with $V(C) = {\lambda}s$.
The trivial
volume estimates show that
the cap $C^{ 1 /(d{\lambda})}$ has volume at most $s$ and is thus
contained in
some set $M(z_i,
15d)$. Then, by Lemma \ref{l:C}, $C$ is contained in $M(z_i,15 d^2
{\lambda})$
which is, by
Lemma \ref{l:B}, a subset of $C^{15 d^2 {\lambda}+1}(z_i) \subset
(C^6(z_i))^{3d^2{\lambda}}=K^{3d^2{\lambda}}_i
(s)$.
\end{pf}

When $P$ is a polytope of volume $V(P)$, the volume of the wet part
$P(v\le s)$ was
determined by Sch\"utt \cite{sch}, as well as by B\'ar\'any and Buchta
\cite{BB}. As $s \to0$,
\[
\frac{V (P(v \leq s V(P)) )}{V(P)}=
\frac{F(P)}{d!d^{d-1}}
s \ln^{d-1} \biggl( \frac1s \biggr) \bigl(1+o(1)\bigr).
\]

Later, we need an estimate for $m(s) $ and $V(P(v \leq s))$,
depending on $P$ only via $F(P)$. Such an estimate follows from
results in B\'ar\'any \cite{Bar1}; see also \cite{Bar7},
formula~(4).
\begin{theorem}\label{th:G} If $P \subset\mathbb{R}^d$ is a polytope with
$V(P)>0$, then
\[
\underline{c} (d) s \ln^{d-1} \biggl( \frac1s \biggr) \leq
\frac{V (P(v \leq s V(P)) )}{V(P)} \leq
\overline{c}(d) F(P) s \ln^{d-1} \biggl( \frac1s \biggr)
\]
and
\[
\underline{c} (d) \ln^{d-1} \biggl( \frac1s \biggr)
\leq m(s V(P)) \leq
\overline{c}(d) F(P)\ln^{d-1} \biggl( \frac1s \biggr)
\]
for $s \leq s_0$, where $\underline{c}(d), \overline{c}(d) >0$ are
constants depending on $d$.
\end{theorem}

The second estimate concerning the number of caps, $m(s)$, follows from
(\ref{ms}).

\section{Plan of proof}\label{sec:plan}

This section explains the basic steps of the proof of Theorem \ref{th:CLTPi}.

\textit{Step} 1. Our proof relies on a precise description of the boundary
of a
convex polytope. The essential ingredients are good bounds on how many
sets $K_i'(s)$ meet a
given cap $C$ of $P$ and on the size of the set visible from $z$ within
$P(v\le T)$. These are obtained in Section \ref{sec:boundary}.

\textit{Step} 2. In what follows, ${\alpha}, {\beta}$ are positive
constants, to be specified later, that depend only on dimension.
Also, we use ``$\operatorname{lln}x$'' as a shorthand for $\ln(\ln
x)$. Define
%
%
\begin{equation}\label{eq:defTseta}
T=T_{\eta}= \frac{{\alpha}\operatorname{lln}{\eta}}\eta
\quad\mbox{and}\quad s=s_{\eta
}=\frac1{\eta\ln^{\beta}\eta} .
\end{equation}
We wish to show that, with high probability, $\Pi_{\eta}$ is sandwiched
between $P(v\ge T)$ and $P(v\ge s)$, that is,
\[
P(v\ge T) \subset\Pi_{\eta} \subset P(v\ge s).
\]
(For technical reasons, we will have to replace $T$ by the slightly
larger $T^*=d6^dT$.) A convenient way to do so is to define a
certain event $A$, which implies sandwiching, and whose complement,
$\overline A$, has very small probability, namely, $\mathbb
{P}(\overline A)
\ll F(P) \ln^{-4d^2}{\eta}$. This will be achieved in Section \ref{sec:sand}.

The basic tool for proving our main result is a central
limit theorem with weakly dependent random variables.
Such an approach has already been used in geometric probability by
Avram and Bertsimas
\cite{AvBe} who also suggested its use in the study of random convex hulls.
For the CLT we are going to use, the weak dependence of random
variables is given by the so-called \textit{dependency graph} which is
defined as follows. Let ${\zeta}_i$, $i \in{\mathcal V}$, be a finite
collection
of random variables. The graph ${\mathcal G}=({\mathcal V},{\mathcal
E})$ is said to be a
dependency graph for ${\zeta}_i$ if, for any pair of disjoint sets
${\mathcal W}_1,
{\mathcal W}_2 \subset{\mathcal V}$ such that no edge in ${\mathcal
E}$ goes between ${\mathcal W}_1$
and ${\mathcal W}_2$, the sets of random variables $\{{\zeta}_i\dvtx i
\in{\mathcal W}_1\}$
and $\{{\zeta}_i\dvtx i \in{\mathcal W}_2\}$ are independent. The
following central
limit theorem with weak dependence is due to Rinott \cite{Ri}. A
slightly weaker version (that would also work here) was earlier proven
by Baldi and Rinott \cite{BalRi}.
\begin{theorem}[(Rinott)]\label{th:rinott}
Let ${\zeta}_i$, $i \in{\mathcal V}$, be random variables having a dependency
graph ${\mathcal G}=({\mathcal V}, {\mathcal E})$. Set ${\zeta}=\sum
_{i \in{\mathcal V}} {\zeta}_i$ and ${\sigma}
^2({\zeta})
= \operatorname{Var}{\zeta}$. Denote the maximal degree of ${\mathcal
G}$ by $D$ and suppose that
$\vert{\zeta}_i - \mathbb{E}{\zeta}_i \vert\leq M$ almost surely.
Then, for every~$x$,
\[
\biggl\vert\mathbb{P}\biggl( \frac{{\zeta}- \mathbb{E}{\zeta}} {\sqrt
{\operatorname{Var}{\zeta}} }\leq x \biggr)
- \Phi(x) \biggr\vert\leq\frac1 {\sqrt{2 \pi}} \frac{DM}{{\sigma}(
{\zeta})} + 16 \frac{\vert{\mathcal V}\vert^{1/2} D^{3/2}
M^2}{{\sigma}^2 (
{\zeta}) } +10 \frac{\vert{\mathcal V}\vert D^2 M^3} {{\sigma}^3
({\zeta}) } .
\]
\end{theorem}

When using this theorem, one has to define the dependency graph and
prove the necessary
properties. Also, we need a lower bound on $\operatorname{Var}{\zeta
}$ (see Theorem \ref
{th:VarPi} below) which
comes from the companion paper \cite{BR}.

\textit{Step} 3. Define a graph ${\mathcal G}$ whose vertex set
${\mathcal V}$ is $\{
1,2,\ldots,m(T)\}$,
where $m(T)$ is the size of the fixed saturated system of points on
$P(v=T)$, as explained
just before the cap covering theorem. The corresponding cap covering
$K_1(T),\ldots,K_{m(T)}(T)$ is indexed by the vertices of ${\mathcal
G}$. Two
vertices $i,j \in{\mathcal V}$
form an edge of ${\mathcal G}$ if the caps $K_i(T)$ and $K_j(T)$ are
``close to
each other,'' in a
well-defined sense. This definition is crucial and will be explained in
Sections \ref{sec:graph} and \ref{sec:CLTunderA}. Also, it will be
shown that the maximal
degree of ${\mathcal G}$ is \mbox{$\ll$}$F(P)^6\operatorname
{lln}^{6(d-1)}{\eta}$.

\textit{Step} 4. Assume that the event $A$ holds which, as mentioned
above, implies
``sandwiching.'' Define the random variables ${\zeta}_i$, $i \in
{\mathcal V}$,
and check that ${\mathcal G}$ is
indeed a dependency graph. The cases of ${\zeta}=V(\Pi_{\eta})$ and
${\zeta}=
f_{\ell}(\Pi_{\eta})$
have to be handled somewhat differently. Next, we check that the
conditions of Rinott's
theorem hold. This will be done in Section \ref{sec:CLTunderA}. This
proves the CLT for ${\zeta}$
given $A$.

\textit{Step} 5. Remove the conditioning on $A$. This is simpler for
${\zeta}=
V(\Pi_{\eta})$, as it is bounded, while ${\zeta}= f_{\ell}(\Pi
_{\eta
})$ is not.
Section \ref{sec:remove} is devoted to this task. The CLT for ${\zeta}$
follows from the CLT for
${\zeta}\vert A$ via the following transference lemma from \cite{bv}, which
has been used in an
implicit form in \cite{Re8} and \cite{vvu2}, and possibly elsewhere.
\begin{lemma} \label{l:trans} Let $\xi_{\eta}$ and $\xi_{\eta}'$
be two
series of random variables
with means $\mu_{\eta}$ and $\mu'_{\eta}$, variances $\sigma_{\eta
}^2$ and
${\sigma}^{\prime2}_{\eta}$, respectively. Assume that there are functions
${\varepsilon}
_1({\eta}), {\varepsilon}_2({\eta}
),{\varepsilon}_3({\eta}), {\varepsilon}_4({\eta}) $, all tending
to zero as ${\eta}$ tends to
infinity, such that:

\begin{longlist}
\item $ |\mu'_{\eta} -\mu_{\eta}| \le{\varepsilon}_1 ({\eta
}) \sigma_{\eta} $;

\item $ | {\sigma}^{\prime2}_{\eta} - \sigma_{\eta} ^2| \le
{\varepsilon}_2({\eta}
) {\sigma}_{\eta} ^2 $;

\item for every $x$, $|\mathbb{P}(\xi'_{\eta} \le x) -
\mathbb{P}(\xi
_{\eta} \le
x)| \le{\varepsilon}_3 ({\eta} ) $;

\item for every $x$,
\[
\biggl| \mathbb{P}\biggl(\frac{\xi'_{\eta}- \mu'_{\eta}}{ \sigma'_{\eta}}
\le x \biggr) -
\Phi
(x) \biggr| \le{\varepsilon}_4({\eta} ) .
\]
\end{longlist}

There is then a positive constant $c$ such that for every $x$,
\[
\biggl| \mathbb{P}\biggl(\frac{\xi_{\eta}- \mu_{\eta}}{ \sigma_{\eta}} \le
x \biggr) -
\Phi(x)
\biggr| \le c\sum_{i=1}^4 {\varepsilon}_i({\eta} ) .
\]
\end{lemma}

The transference lemma asserts that if $\xi'_{\eta}$ satisfies the CLT
(the fourth condition)
and $\xi_{\eta}$ is sufficiently close to $\xi'_{\eta}$ in distribution
(the first three
conditions), then $\xi_{\eta}$ also satisfies the CLT.
\begin{Remark*}
In \cite{bv}, the transference lemma is stated with
$\sigma'_{\eta}$ and
${\sigma}^{\prime2}_{\eta}$ on the right-hand side of conditions (i) and (ii).
It is easy to see
that the present conditions imply those involving $\sigma'_{\eta}$:
(ii) shows that
${\sigma}^{\prime2}_{\eta}/ \sigma_{\eta}^2$ tends to $1$ as $n \to\infty$.
Thus, $\sigma_{\eta}^2 <
2{\sigma}^{\prime2}_{\eta}$ for large enough $n$. Then (ii) implies
$|{\sigma}^{\prime2}_{\eta}-\sigma_{\eta}^2| \le2{\varepsilon}_2({\eta}
) {\sigma}^{\prime2}_{\eta}$ and, similarly,
(i) implies $|\mu'_{\eta}-\mu_{\eta}| \le\sqrt{2}{\varepsilon
}_1({\eta} ) \sigma
'_{\eta}$.
\end{Remark*}

To apply the central limit theorem and the transference lemma, we need
a lower bound on $\operatorname{Var}
{\zeta}$. In the companion paper \cite{BR}, we prove a lower bound for
general convex bodies in
terms of the volume of the floating body: Theorem 3.1 in \cite{BR} says
that the variance of
$V(\Pi_\eta)$ is bounded from below by $\eta^{-1} V(K(v \leq\eta
^{-1}))$ and $\operatorname{Var}
f_{\ell}(\Pi_\eta)$ is bounded by $ \eta V(K(v \leq\eta^{-1}))$. Using
Theorem \ref{th:G},
this gives the following result.
\begin{theorem}\label{th:VarPi} Assume that $P$ is a polytope of volume
$1$. Then
\begin{eqnarray*}
F(P)\eta^{-2} \ln^{d-1} \eta &\ll& \operatorname{Var}V(\Pi_{\eta}),
\\
F(P) \ln^{d-1} \eta &\ll& \operatorname{Var}f_{\ell}(\Pi_{\eta}) .
\end{eqnarray*}
\end{theorem}

Here, we use Vinogradov's ``$\gg$'' notation, that is, we write
$f({\eta
} ) \gg
g({\eta} )$ if there is a constant $c>0$, independent of $\eta$, such that
$cf({\eta} )>|g({\eta} )|$ for all $\eta\geq\eta_0 $. The constants
$c$ and $\eta_0$
may, and usually do, depend on the dimension, but not on
$K$.

The main achievements of this paper, besides the central limit
theorems, are the
precise sandwiching of $\Pi_{\eta}$, the novel definition of the
dependency graph and the
proof that its maximal degree is bounded by a power of $\operatorname
{lln}{\eta}$.
The latter is based on
structural properties of the wet part $P(v \le t)$ for polytopes.

\section{On the boundary structure of convex polytopes}\label{sec:boundary}

In this section, we state some facts about the boundary structure of
the polytope $P$ and its
floating body. All proofs in this section, except for those of
Claim \ref{cl:simp} and
Lemma \ref{l:visible}, which are given here, are postponed to Section
\ref{sec:aux}.

So, the polytope $P$ is fixed and its volume is $1$. We need to
consider two parameters, $T$ and $s$, which have already been
defined in (\ref{eq:defTseta}). However, this is not important for the time
being; we only assume that $2s\le T$, say.

Let $z \in P$ be a point with $v(z) \le T$ and write $[x,z]$ for the
closed segment joining $z$ and a point $x$. The following definition is
crucial, also having been used
by Vu \cite{vvu1}. Set
\[
S(z,T)=\{ x \in P\dvtx[x,z] \cap P(v \ge T) = \varnothing\}.
\]
This is the set of points that are \textit{visible from} $z$ \textit{within}
$P(v < T)$. We are
interested in the size of $S(z,T)$.

We again use the notation $g(s) \ll f(s)$ if $|g(s)| < cf(s)$ for
all $0< s \leq t_0 $ with constants $c$ and $t_0$ depending on the
dimension, but not on the underlying convex set.
\begin{lemma} \label{l:volS}
If $0 < v(z) \leq\frac12 $, $2 v(z) \leq T$, then
\[
V(S(z,T)) \ll F(P) T \ln^{d-1} \biggl( \frac T {v(z)} \biggr).
\]
\end{lemma}

Note that since $S(z,T) \subset P(v\le T)$, Theorem \ref{th:G}
immediately implies the inequality $V(S(z,T)) \ll F(P) T
\ln^{d-1} (1/T)$. The improvement from $1/T$ to $T/v(z)$ is
significant in the range we are interested in.

Consider the economic cap covering from Theorem \ref{th:capcov} for
$P(v\le s)$, $s \leq s_0$, where $s_0$ is defined in (\ref{def:s0}).
The caps $K_i(s)$ come from a saturated system
$Z(s)=\{z_1,\ldots,z_{m(s)}\} \subset P(v=s)$ which is fixed together
with $P$ and $s$, as agreed just before the cap covering theorem was presented.
We want to know how many $z_i \in Z(s)$ can be contained in a fixed
cap $C$ of volume $T$.
\begin{lemma} \label{l:zinC}
Assume that $C$ is a cap of $P$ of volume $T$. Then, for $0 < s \leq
s_0$, $ 2s \leq T$, we have
\[
\vert Z(s) \cap C\vert\ll F(P) \ln^{d-1} \biggl(\frac Ts \biggr).
\]
\end{lemma}

Next, consider the economic cap covering theorem for $P(v\le T)$. The
saturated system $Y(T)=\{y_1,\ldots,y_{m(T)}\}$ on $P(v=T)$ is again
fixed and so are the corresponding covering caps $K_j(T)$. [We use
the notation $Y(T)$, $y_j(T)$ and $m(T)$ in order to avoid confusion
with $Z(s),z_i(s)$ and $m(s)$.] We will need a bound on the number
of those $y_j \in Y(T)$ for which $K_j^{\lambda}(T)$ contains a fixed $z
\in P(v=s)$. Here, ${\lambda}$ is a constant that depends only on $d$.
\begin{lemma} \label{l:zinK}
Let ${\lambda}\geq1$ be a constant depending only on $d$. Assume that
$0\le
2s \le T \le(6 {\lambda})^{-d} s_0$. If $z \in P(v=s)$, then
\[
| \{ y_j \in Y(T)\dvtx z \in K_j^{\lambda}(T) \} |
\ll F(P) \ln^{d-1} \biggl(\frac Ts \biggr) .
\]
\end{lemma}

The constant in $\ll$ depends on ${\lambda}$ and, thus, again, only
on the dimension.

We will also need a bound on the number of points $z_j \in Z(s)$
that are contained in $S(z,T)$ when $z \in P(v=s)$.
\begin{lemma} \label{l:nkl} Assume that $z \in P(v=s)$ and $0 < s \leq
s_0, 2s \leq T $.
Then
\[
\vert Z(s) \cap S(z,T) \vert\ll F(P) \ln^{d-1} \biggl(\frac T{s} \biggr).
\]
\end{lemma}

The following fact will be needed in the sandwiching step and
concerns convex hulls of random points in $K'_i(T)$, the small sets
in the cap covering theorem. Set $T^*=d6^dT$, $T \leq s_0$. In each
$K'_i(T)$, choose a point $x_i$ arbitrarily.
\begin{claim}\label{cl:convxi1} Under the above conditions,
\[
P(v \geq T^*) \subset\bigl[x_1, \ldots, x_{m(T)} \bigr] .
\]
\end{claim}

We mention in passing that the caps $K^{\gamma}_i (T)$ cover $P(v\le
T^*)$, where ${\gamma}=3d^36^d$:
%
%
\begin{equation}\label{eq:T*}
P(v\le T^*) \subset\bigcup_1^{m(T)} K_i^{\gamma}(T).
\end{equation}
This follows directly from Claim \ref{cl:lambda}.

The system $Z(s)=\{z_1,\ldots,z_{m(s)}\}$ on $P(v=s)$ is saturated,
so, for each $a \in P(v=s)$, there is a $z_i$ with $M(z_i,\frac12)\cap
M(a, \frac12)\ne\varnothing$. For each $a$, we fix such a $z_i$ and denote
it by $z(a)$.
\begin{claim}\label{cl:simp} If a cap $C$ contains the point $a \in P(v=s)$,
then $M(z(a),1) \subset C^6$.
\end{claim}
\begin{pf}
This is very simple. As $z(a)$ satisfies
$M(z(a), \frac12)\cap M(a, \frac12) \ne\varnothing$ by definition,
it follows that Lemmas \ref{l:A} and \ref{l:B} imply that
\[
M(z(a),1)=P\cap M(z(a),1)\subset P\cap M(a,5) \subset C^6.
\]
\upqed\end{pf}

The following lemma helps to bound the maximal degree of the
dependency graph.
\begin{lemma} \label{l:visible} Assume that $a,b\in P(v=s)$ and the
segment $[a,b]$ is disjoint from $P(v\ge T)$. The segment
$[z(a),z(b)]$ is then disjoint from $P(v\ge T^*)$.
\end{lemma}
\begin{pf}
Both $[a,b]$ and $P(v\ge T)$ are convex, so they can be
separated by a hyperplane since they are disjoint. This hyperplane
cuts off a cap, say $C$, from $K$ containing $[a,b]$ and disjoint
from $P(v\ge T)$. So, $V(C) \leq dT$, by Lemma \ref{l:E}. Further,
Claim \ref{cl:simp} implies that $z(a),z(b) \in C^6$. Consequently,
$[z(a),z(b)] \subset C^6$ and $V(C^6) \le d6^dT=T^*$ follows from
the trivial volume estimate.
\end{pf}

\section{Sandwiching $\Pi_{\eta}$}\label{sec:sand}

Recall that the Poisson polytope, $\Pi_{\eta}$, is the convex hull of
$X \cap P$, where $X=X(\eta)$ is a Poisson point process of intensity
$\eta$. We
are going to use the well-known fact that, with high probability, the
boundary of $\Pi_{\eta}$ is contained in a small strip close to the
boundary of $P$. Results of this type have been proven in \cite{BD}
and \cite{vvu1}. Here, we need a slightly different, perhaps more
refined, estimate.

We make (\ref{eq:defTseta}) more precise and set
\[
T=T_{\eta}= {\alpha}\frac{ \operatorname{lln}{\eta}} \eta
\qquad\mbox{with }{\alpha}
=(6d)^d(4d^2+d-1).
\]
In the following, we assume that $\eta\geq\eta_0$, where $\eta_0$ is
chosen such that $T \leq s_0$, with $s_0$ defined as in (\ref{def:s0}).
Let $Y(T)$ be the fixed saturated point set $\{y_1,\ldots,$
$y_{m(T)}\}$ on $P(v= T )$ according to Theorem \ref{th:capcov}. We
get an economic cap covering with caps $K_j(T)$ and half Macbeath
regions $K'_j(T)$, $j=1, \ldots,m(T)$. To simplify notation,
set $K_j=K_j(T)$, $K_j'=K_j'(T)$ and $m_{\eta}=m(T)$.

Let $A'$ be the event that each $K'_j$ contains at least one point
of $X$, the Poisson point process with intensity $\eta$. Since the
number of points in $K'_j$ is Poisson distributed with parameter $\eta
V(K'_j)$, from the fact that $(6d)^{-d} T \le V(K'_j) \le2^{-d}T$, we have
\[
\mathbb{P}(K'_j \cap X=\varnothing)= e^{- \eta V(K'_j)} \leq e^{-(6d)^{-d}
\eta
T }.
\]
Let $\overline{A}{}'$ denote the complement of the event $A'$. By
Theorem \ref{th:G}, $m_{\eta} \ll F(P)\times \ln^{d-1} \eta$, so, by
Boole's inequality,
%
%
\begin{equation}\label{A'pi}\quad
\mathbb{P}(\overline{A}{}')\leq m_{\eta} e^{-(6d)^{-d} \eta T} \ll F(P)
(\ln\eta)^{-(6d)^{-d}{\alpha}+ d-1} = F(P) \ln^{-4d^2} \eta
\end{equation}
follows from the choice of ${\alpha}$.

For later reference, we note that
%
%
\begin{equation}\label{eq:A'pi>}
\mathbb{P}(K'_j \cap X=\varnothing)\geq e^{-2^{-d} \eta T }
= \ln^{-2^{-d}{\alpha}} \eta\geq\ln^ {-(3d)^ {d+2}} \eta.
\end{equation}

Now, Claim \ref{cl:convxi1} and (\ref{A'pi}) show that,
with high probability, $\Pi_{\eta}$ contains the floating body
$P(v \geq T^*)$. (Recall that $T^*= d 6^d T$.)
\[
\mathbb{P}\bigl(\Pi_{\eta} \mbox{ does not contain }P(v \geq T^*)\bigr)\le
\mathbb{P}
(\overline{A}{}') \ll F(P) \ln^{-4d^2} \eta.
\]

This is the first half of the sandwiching. For the second half, we make
the definition of $s_n$ in (\ref{eq:defTseta}) more precise and set
\[
s= s_{\eta}= \frac1 { \eta\ln^{\beta} \eta}\qquad \mbox{where
}{\beta}=4d^2+d-1.
\]
We claim that, with high probability, $P(v \leq s)$ contains no point
of $X$. Indeed, $ \eta V(P(v \leq s)) \ll F(P) (\ln{\eta})^{-{\beta}
+d-1}$, by
Theorem \ref{th:G}, and we get
%
%
\begin{equation}\label{eq:empty}
\mathbb{P}\bigl(X\cap P(v \le s)\ne\varnothing\bigr)=1- e^{- \eta V(P(v \leq
s))} \ll
F(P) \ln^{-4d^2} \eta.
\end{equation}

We have just proven that $\Pi_{\eta}$ is sandwiched between $P(v\ge
s)$ and $P(v\ge T^*)$ with high probability:
\[
1-\mathbb{P}\bigl( P(v\ge T^*)\subset\Pi_{\eta} \subset P(v\ge s) \bigr)
\ll F(P) \ln^{-4d^2} \eta.
\]

\begin{figure}

\includegraphics{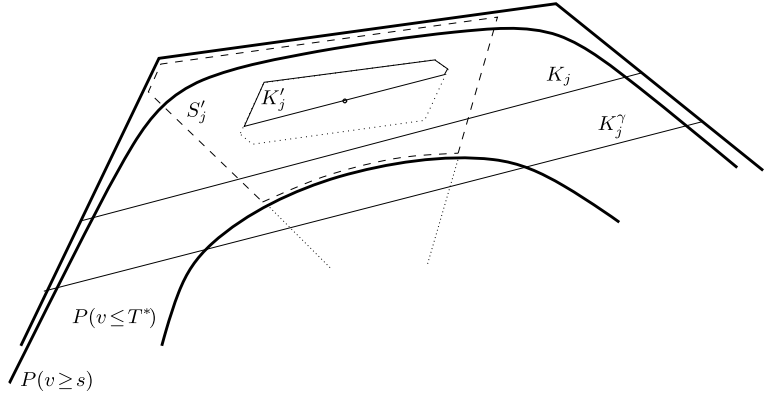}

\caption{Definition of $S_j'$.}
\label{figure1}
\end{figure}

The proof of the CLT for $V(\Pi_{\eta})$ could be achieved via conditioning on
$A'$. For
$f_{\ell}(\Pi_{\eta})$, we need a stronger condition, to be called $A$,
which will also work for~$V(\Pi_{\eta})$. Set
\[
{\gamma}= 3 d^3 6^d.
\]
For $j=1, \ldots,m_{\eta}$, let $S_j =S_j(T)$ be pairwise internally
disjoint closed sets with $\bigcup S_j=P$,
$K'_j \subset S_j$ and $S_j \cap P(v \leq T^* ) \subset K^{\gamma}_j $.
[Recall, from Claim \ref{cl:lambda}, that the sets
$K^{\gamma}_j$ cover $P(v \leq T^*)$.] Set $S'_j= S'_j(T)=S_j \cap P(v
\leq T^* )$, see Figure \ref{figure1}.

%
%

Before defining $A$, observe that the expected number of points of $X$
lying in $S_j'$ is $ \eta V(S'_j)$. Trivial volume
estimates show that
%
%
\begin{equation}\label{eq:Sj}
(6d)^{-d} {\alpha}\operatorname{lln}{\eta}
\le\eta V(K'_j) \le\eta V(S'_j) \le\eta V(K_j^\gamma) \le(6\gamma
)^d {\alpha}
\operatorname{lln}{\eta}.
\end{equation}

Define $A$ to be the event that each $K'_j$ contains at least one
point, $P(v \leq s)$ contains no point and each $S'_j$ contains
at most $3(6{\gamma})^d {\alpha}\operatorname{lln}{\eta}$ points of
$X$ ($j=1, \ldots,
m_{\eta}$).
The following two claims are essential for our proof. We collect
the properties of $\Pi_{\eta}$ given the event $A$ and estimate the
probability of $A$.
\begin{claim}\label{cl:Aimpl}
Given $A$, we have $P(v\ge T^*)\subset\Pi_{\eta} \subset P(v\ge s) $ and
$| P(v\le T^*) \cap X |\ll F(P) \ln^ {d-1} \eta
\operatorname{lln}{\eta
} $.
\end{claim}
\begin{pf}
This follows immediately from the definition of $A$ and
from the
estimate on the volume of $P(v\le T^*)$.
\end{pf}
\begin{claim}\label{cl:PA}
$ \ln^ {-(3d)^ {d+2}} \eta\ll\mathbb{P}( \overline{A} ) \ll F(P)
\ln
^{-4d^2} \eta$.
\end{claim}
\begin{pf}
The lower bound follows from $\mathbb{P}(\overline
{A})\ge\mathbb{P}
(K'_1(T)\cap
X=\varnothing)$ and from (\ref{eq:A'pi>}). For the upper bound,\vspace*{-1pt} recall
(\ref{A'pi}) and
(\ref{eq:empty}): in (\ref{A'pi}) we showed that $K'_j \cap
X=\varnothing
$ for some $j$ has
probability \mbox{$\ll$}$ F(P) \ln^{-4d^2} \eta$; inequality (\ref{eq:empty})
shows that $X \cap
P(v\le s)\ne\varnothing$ with probability \mbox{$\ll$}$F(P) \ln^{-4d^2}
\eta$.

So, we only have to estimate the probability that for some $j$, the set
$S'_j$ contains more
than $3(6{\gamma})^d {\alpha}\operatorname{lln}{\eta} \geq\eta
V(S'_j)$ points. Let $N$ denote
a Poisson random
variable with parameter $p$. Then (see, e.g., \cite{Re8})
\[
\mathbb{P}( N \geq3 p) \leq\frac3 {3-e} e^{- p}.
\]
This inequality implies, by setting $p= \eta V(S'_j)$, that the
probability that $S'_j$
contains more than $3(6{\gamma})^d {\alpha}\operatorname{lln}{\eta}
\geq3p$ points from $X$ is
bounded from above by
\[
\frac3 {3-e} e^{-p} \le\frac3 {3-e} \exp(-
(6d)^{-d} {\alpha}\operatorname{lln}{\eta} ) \ll\ln^{-(4d^2 +d
-1)} \eta.
\]
Combining this with the bound $m_{\eta} \ll F(P)\ln^{d-1} \eta$ from
Theorem \ref{th:G}
completes the proof.
\end{pf}

Set
\[
U= U_{\eta}=\frac{\ln{\eta}} \eta\quad\mbox{and}\quad U^*= d6^d U.
\]
Since we are assuming that $V(P)=1$, Theorem \ref{th:G} tells us that
\[
b_1 \frac{\ln^d \eta} \eta\le V\bigl(P(v \le U ^*)\bigr) \le b_2 F(P) \frac
{\ln^d \eta} \eta
\]
with positive constants $b_1,b_2$ depending only on $d$. Let $B$ be
the event that $P(v \ge U^*) \subset\Pi_{\eta}$ and that $P(v \le U^*)$
contains at most $3b_2F(P) \ln^d \eta$ points from $X$. The
following estimate will be useful in Section \ref{sec:remove}. Its
proof is similar to those above, actually even simpler (as there is no need
to worry about ${\alpha}$), and is therefore left to the reader.
\begin{lemma}\label{l:Bpi}
$\mathbb{P}(\overline{B}) \ll F(P) \eta^{-3d}$.
\end{lemma}
%

\section{The dependency graph}\label{sec:graph}

It is high time to define the dependency graph ${\mathcal G}=
({\mathcal V}, {\mathcal E})$.
The values of $s,T$ and $T^*$ have been given in the previous section. The
sets $K_i=K_i(T)$ and $K'_i=K'_i(T)$ come from the cap covering
theorem. The vertex set, ${\mathcal V}$, of the dependency graph is just
$\{1, \ldots, m_{\eta}\}$.

Define the set $L_i$ as the union of all $S'_k$ such that there are points
%
%
\begin{eqnarray}\label{def:Li}
a \in S'_i\cap P(v\ge s),\qquad b \in S'_k\cap P(v\ge
s)\nonumber\\[-8pt]\\[-8pt]
\eqntext{\mbox{with } [a,b] \mbox{ disjoint from }P(v\ge T^*).}
\end{eqnarray}
Note that $S_i'\subset L_i$ for all $i$. Also, $S_k'\subset L_i$
holds if and only if $S_i'\subset L_k$. Now, distinct vertices $i,j
\in{\mathcal V}$ form an edge in ${\mathcal G}$ if $L_i $ and $ L_j $
contain at
least one set $S'_k$ in common,
%
%
\begin{equation}\label{def:depgr}
ij \in{\mathcal E}\quad\Leftrightarrow\quad\exists k \in\{ 1, \ldots, m_{\eta
} \}\qquad\mbox{such that }
S'_k \subset L_i \cap L_j.
\end{equation}

That this defines a dependency graph for the suitably chosen
random variables is proved later, in Lemma \ref{le:rinott}. The
main result of this section is an upper bound on the maximal
degree $D$ in ${\mathcal G}$.

Some preparation is needed. We need a bound on the number of sets $S'_k
\subset L_i$.
\begin{lemma}\label{l:SkLi}
$ | \{ k \dvtx S'_k \subset L_i \} | \ll F(P)^3(\operatorname{lln}\eta
)^{3(d-1)}$.
\end{lemma}
\begin{pf}
We first show that if $S'_k \subset L_i$, then there are also points
%
%
\begin{eqnarray}\label{eq:ab'}
a' \in K_i^\gamma\cap P(v=s) ,\qquad b' \in K_k^\gamma\cap
P(v=s) \nonumber\\[-8pt]\\[-8pt]
\eqntext{\mbox{with } [a',b'] \mbox{ disjoint from
}P(v\ge T^*).}
\end{eqnarray}
To simplify notation, we write $C'=C \cap P(v\ge s)$ when $C$ is a cap
of $P$. Clearly, $C'$
is a cap of $P(v\ge s)$. We are going to use the fact that if two caps
of a convex body have
a point in common, then they also have a point in common from the
boundary of the convex set.

Since the segment $[a,b]$ is disjoint from $P(v\ge T^*)$, there is a
cap $C$, also disjoint from $P(v\ge
T^*)$, such that $[a,b] \subset C'$. Now,
\[
a \in C' \cap {K{}^{\gamma}_i}{}',\qquad b \in C' \cap K_k^{{\gamma}}{'}.
\]
Every one of the two sets above is a nonempty intersection of two caps
of $P(v\ge s)$. So,
each has a point, $a'$ and $b'$, respectively, on the boundary of
$P(v\ge s)$, which is
$P(v=s)$. As the segment $[a',b'] \subset C$, it is disjoint from
$P(v\ge T^*)$, which proves
(\ref{eq:ab'}).

Recall that a saturated system $Z(s)=\{z_1,\ldots,z_{m(s)}\}$ has been
chosen in $P(v=s)$.
Also, for each $x \in P(v=s)$, we fixed a point $z(x) \in Z(s)$ so that
$M(x,{ \frac12})\cap
M(z(x),{ \frac12})\ne\varnothing$. We have points $a',b' \in P(v=s)$
satisfying
(\ref{eq:ab'}). Claim \ref{cl:simp} shows the existence of points
$z(a'),z(b')\in Z(s)$ such
that $z(a') \in K_i^{6{\gamma}}$, $z(b') \in K_j^{6{\gamma}}$ and, by
Lemma \ref
{l:visible}, the segment
$[z(a'),z(b')]$ is disjoint from $P(v\ge T^o)$, where $T^o=d6^dT^*$.

We bound the number of sets $S'_k$ in $L_i$ in three steps.
In view of Lemma
\ref{l:zinC}, with $C=K_i^{6{\gamma}}$, we have
\[
\vert Z(s) \cap K_i^{6{\gamma}}\vert
\ll F(P) \ln^{d-1} \biggl(\frac{V(K_i^{6{\gamma}})}s \biggr) \ll F(P)
\operatorname{lln}^{d-1} \eta,
\]
where the upper bound for $V(K_i^{6{\gamma}})$ comes from (\ref
{eq:Sj}). This
is an upper bound on
how many $z(a') \in K^{6{\gamma}}_i$ there can be, given that the segment
$[z(a'), z(b')]$ starts
at $K_i^{6{\gamma}}$.

In the second step, we estimate, for a fixed $z(a')$, the number of
$z(b') \in Z(s)$ such that $[z(a'),z(b')]$ is disjoint from $P(v\ge T^o)$.
All such $z(b')$ lie in $S(z(a'),T^o)$. So, by Lemma \ref{l:nkl}, the
number of such $z(b')$ is
\[
\mbox{$\ll$}F(P) \ln^{d-1} \biggl(\frac{T^o}s \biggr) \ll F(P) \operatorname{lln}^{d-1}
\eta.
\]

In the third step, we estimate the number of $K_j^{6{\gamma}}$ that
contain a fixed $z(b') \in Z(s)$. Lemma \ref{l:zinK} implies, with
${\lambda}=6{\gamma}$, that this number is
\[
\mbox{$\ll$} F(P)\ln^{d-1} \biggl( \frac T s \biggr) \ll F(P) \operatorname{lln}^{d-1}
\eta.
\]

This argument shows that for a set $S'_i$, there are at most \mbox{$\ll$}$
F(P)^3 \operatorname{lln}^{3(d-1)} \eta$ sets $S_k'$ which can be
connected by
some segment $[a', b']$. Since every set $S_k' \subset L_i$ is
connected to $S'_i$ by some segment, the number of sets $S'_k$ in $L_i$ is
\mbox{$\ll$}$ F(P)^3 \operatorname{lln}^{3(d-1)} \eta$.
\end{pf}

The following result gives the upper bound on the maximal degree $D$.
\begin{theorem}\label{th:degree} $D \ll F(P)^6(\operatorname
{lln}{\eta})^{6(d-1)}$.
\end{theorem}
\begin{pf}
By (\ref{def:depgr}), we have $ij \in{\mathcal E}$
if $L_i
\cap L_j $ contains some set $S'_k$. Clearly, if $S'_k \subset
L_j$, then, by the definition (\ref{def:Li}), we also have $S'_j
\subset
L_k$. Thus, $ij \in{\mathcal E}$ if there is some $k$ such that $S'_k
\subset L_i, S'_j \subset L_k$, which gives
\[
D \leq{\max_{i} \sum_{ k \dvtx S'_k \subset L_i } }| \{ j \dvtx S'_j
\subset
L_k \} | .
\]
Combined with Lemma \ref{l:SkLi}, this gives the bound on the
degree of ${\mathcal G}$.
\end{pf}

Thus, the graph ${\mathcal G}$ has been defined and its maximal degree
has been
bounded. In the next section, we define the random variables ${\zeta}_i$
and show that ${\mathcal G}$ is a dependency graph.

\section{The central limit theorem under condition $A$}\label{sec:CLTunderA}

\mbox{}

\begin{pf*}{Proof of the CLT for $V(\Pi_{\eta})|A$}
We introduce $m_{\eta}$ random
variables ${\zeta}_j$ in the following way. For simpler notation, we
keep writing $K'_j$ for $K'_j(T)$, $K_j$ for $K_j(T)$, $S_j$ for $S_j(T)$ and
$S'_j$ for $S'_j(T)$.
We define ${\zeta}_j$ as the missed volume in the set $S_j$,
\[
{\zeta}_j = V(S_j) - V(S_j \cap\Pi_{\eta}) ,
\]
and ${\zeta}$ as the missed volume in the polytope $P$,
\[
{\zeta}= \sum_{j=1}^{m_{\eta}} {\zeta}_j =V(P)-V(\Pi_{\eta
})=1-V(\Pi_{\eta}).
\]

In order to prove the CLT for $V(\Pi_{\eta})|A$, we simply check the
conditions of Rinott's theorem. We start with the weak independence
condition.
\begin{lemma}\label{le:rinott}
Given disjoint subsets ${\mathcal W}_1,{\mathcal W}_2$ of ${\mathcal
V}$ with no edge between
them, the random variables $\{{\zeta}_i\dvtx i \in{\mathcal W}_1\}$
are independent
of the random variables $\{{\zeta}_j\dvtx j \in{\mathcal W}_2\}$
under the
conditional distribution of $X$ given that $A$ holds.
\end{lemma}
\begin{pf}
Under condition $A$, the boundary of $\Pi_{\eta}$
lies in
$P(s<v\le
T^*)$ and, thus, ${\zeta}_j = V(S'_j) - V(S'_j \cap\Pi_{\eta}) $. The
intersection $S'_i \cap
\Pi_{\eta}$ is determined by the facets [$(d-1)$-dimensional faces] of
$\Pi_{\eta}$
intersecting $S'_i$. These facets are determined by their vertices.
Thus, all vertices that
may determine a facet that intersects $S'_i$ are contained in $L_i$. In
other words, $ S'_i
\cap\Pi_{\eta}$ is the same as the intersection of $ S'_i$ with the
convex hull of $ X \cap
L_i$.

Now, set $L^k=\bigcup_{i \in{\mathcal W}_k}L_i$ for $k=1,2$. By definition,
$L^1$ and $L^2$ are unions
of sets $S'_k$ and have disjoint interiors.
Given $A$, the ${\zeta}_i, i \in{\mathcal W}_1$, are determined
by $L^1 \cap X$ and the ${\zeta}_i, i \in{\mathcal W}_2$, are
determined by $L^2
\cap X$.
Since $L^1 \cap X$ is independent of $L^2 \cap X$, conditional on $A$
and otherwise, the claim follows.
\end{pf}

We have to check two more conditions of Rinott's theorem.
\begin{claim}\label{cl:Mvol} Under condition $A$, $M= \max\| {\zeta
}_j \|
_{\infty} \ll(\operatorname{lln}{\eta})/ \eta$.
\end{claim}
\begin{pf}
This is very simple: ${\zeta}_j \leq V(S'_j) \ll T
\ll
(\operatorname{lln}
{\eta})/{ \eta}$.
\end{pf}
\begin{claim}\label{cl:var*}
For $\ln\eta\gg F(P)^{1/d^2}$, we have $\operatorname{Var}(V(\Pi
_{\eta})|A) \gg
F(P)\eta^{-2} \times\break\ln^{d-1}
\eta$.
\end{claim}

This claim is an easy corollary of Theorem \ref{th:VarPi} and (\ref
{eq:varVA}) from the next
section.

Bounds on $|{\mathcal V}|, D, {\zeta}_j$ and $\operatorname
{Var}{\zeta}=\operatorname{Var}(V(\Pi_{\eta})|A)$ have been
established. Rinott's theorem can be applied.
For $\ln\eta\gg F(P)^{1/d^2}$, the dominating error term is
\[
\frac{|{\mathcal V}|D^2M^3}{\operatorname{Var}(V(\Pi_{\eta})|A)^{
3/2 }} \ll
F(P)^{11.5}\frac
{(\operatorname{lln}{\eta})^{12d-9}}{(\ln{\eta})^{({d-1})/2}}
\]
as a simple computation shows. If $\ln\eta$ equals $ F(P)^{1/d^2}$,
then the right-hand side is
already $\gg$1, which proves that this error term is valid for all
$\eta$.
\end{pf*}
\begin{pf*}{Proof of the CLT for $f_{\ell}(\Pi_{\eta})|A$} The dependency
graph remains the same.
The random variables ${\zeta}_i$ are to be defined, just as in \cite{Re8},
in the following way.
Let $F$ be an $\ell$-dimensional face of $\Pi_{\eta}$ having
$f_0(S_i,F)$ vertices in $S_i$
and set
\[
{\zeta}_i= \frac1{\ell+1} \sum_{\mathrm{all}\ F}f_0(S_i,F).
\]
Since, with probability one, no point from $X$ lies in two $S_j$, and
each face $F$ is a
simplex with probability one, the sum of the ${\zeta}_i$ is equal to
$f_{\ell
}(\Pi_{\eta})$ almost
surely. The analog of Lemma \ref{le:rinott} for the new variables
${\zeta}
_i$ is proved in the same
way.

We need to bound $\max\| {\zeta}_i \|_{\infty}$ from above and,
also, $\operatorname{Var}
{\zeta}
=\operatorname{Var}
f_{\ell}(\Pi_{\eta})$ from below.
\begin{claim}\label{cl:varface}
For
$\ln\eta\gg F(P)^{1/d}$,
we have
$\operatorname{Var}(f_{\ell}(\Pi_{\eta})|A) \gg F(P) \ln^{d-1}
\eta$.
\end{claim}

Again, this follows from Theorem \ref{th:VarPi} and (\ref{eq:varfA})
in the next section.
\begin{claim}\label{cl:Mface}$M=\max\| {\zeta}_i \|_\infty\ll
F(P)^{3d}(\operatorname{lln}{\eta})^{3d^2}$.
\end{claim}
\begin{pf}
(Similar to the one in Reitzner \cite{Re8}.) Condition
$A$ ensures that all vertices of $\Pi_{\eta}$ lie in $P(s< v \le T^*)$.
As we have seen in the proof of Lem\-ma~\ref{le:rinott}, each face $F$
intersecting $S'_i$ has all of its vertices in $L_i$: if $x\in S'_i$
and $y\in S'_j$ are vertices of $F$, then $y \in L_i$.
Under condition $A$, $S'_j$ contains \mbox{$\ll$}$\operatorname{lln}{\eta}$ points
from~$X$. Thus, the number of vertices contributing to ${\zeta}_i$ is
\mbox{$\ll$}$ F(P)^3(\operatorname{lln}{\eta})^{(3d-2)} $, by Lemma \ref{l:SkLi}.

The number of $\ell$-faces (actually, all subsets of
size $\ell+1$) on this many vertices is \mbox{$\ll$}$(F(P)^3(\operatorname
{lln}{\eta
})^{(3d-2)} )^{\ell+1}$.
Each such $\ell$-face contributes at most $1$ to the value of
${\zeta}_i$. Consequently,
\[
{\zeta}_i \ll\bigl(F(P)^3(\operatorname{lln}{\eta})^{(3d-2)} \bigr)^{\ell+1}
\ll F(P)^{3d}(\operatorname{lln}{\eta})^{d(3d-2)}
\]
since $\ell+1 \le d$.
\end{pf}

All conditions of Rinott's theorem have been established. The
dominating error term is again the third one and we get the CLT for
$f_{\ell}(\Pi_{\eta})\vert A$ with error term
\[
\frac{|{\mathcal V}|D^2M^3}{\operatorname{Var}(f_{\ell}(\Pi_{\eta
})|A)^{3/2}} \ll
F(P)^{15d}\frac{(\operatorname{lln}{\eta})^{15d^2}}{(\ln{\eta
})^{({d-1})/2}}
\]
as a simple computation shows.
\end{pf*}

\section{Removing the conditioning}\label{sec:remove}

We are going to use the transference Lem\-ma~\ref{l:trans}.
\begin{lemma}\label{l:transV} The random variables $\xi_{\eta}=V(\Pi
_{\eta})$ and
$\xi'_{\eta}=V(\Pi_{\eta})|A$ satisfy the conditions of Lemma \ref{l:trans}
with
\[
\sum{\varepsilon}_i ({\eta} ) \ll F(P)^ {11.5} \ln^{-(
{d-1})/{2}+o(1)} \eta.
\]
\end{lemma}
\begin{lemma}\label{l:transfs} The random variables $\xi_{\eta
}=f_{\ell
}(\Pi_{\eta})$
and $\xi'_{\eta}=f_{\ell}(\Pi_{\eta})|A$ satisfy the conditions of
Lemma \ref{l:trans} with
\[
\sum{\varepsilon}_i ({\eta} ) \ll F(P)^ {15d} \ln^{-(
{d-1})/{2}+o(1)} \eta.
\]
\end{lemma}

In both cases, the fourth condition of the transference lemma has been
proven in the previous section with ${\varepsilon}_4 \ll
F(P)^{11.5} \ln^{-(d-1)/2 +o(1)} \eta$ for the case of volume and with
${\varepsilon}_4 \ll F(P)^{15d}
\ln^{-(d-1)/2+o(1)} \eta$ for the number of faces. So, our main theorem
for $\Pi_{\eta}$
follows once the first three conditions of the transference lemma have
been checked for the
volume and for the number of faces. We will make use of the following
simple claim.
\begin{claim}\label{cl:simple} If $\zeta$ is a nonnegative
random variable and $A$ is an event, then
\[
|\mathbb{E}(\zeta) -\mathbb{E}(\zeta|A)| \le\bigl(\mathbb{E}(\zeta
|A)+\mathbb{E}(\zeta|\overline
{A}) \bigr)\mathbb{P}(\overline{A}).
\]
\end{claim}
\begin{pf}
It is clear that $\mathbb{E}(\zeta)=\mathbb
{E}(\zeta|A)\mathbb{P}(A)
+\mathbb{E}(\zeta|\overline{A})\mathbb{P}(\overline{A})$.
Replacing $\mathbb{P}(A)$ by
$1-\mathbb{P}(\overline{A})$ here gives
\[
\mathbb{E}(\zeta) -\mathbb{E}(\zeta|A) = \bigl(-\mathbb{E}(\zeta
|A)+\mathbb{E}(\zeta|\overline
{A}) \bigr)\mathbb{P}(\overline{A})
\]
and the claim follows.
\end{pf}
\begin{pf*}{Proof of Lemma \ref{l:transV}}
We need some preparations. We
use Claim \ref{cl:simple} with $\zeta=1-V(\Pi_{\eta})$. We first estimate
$\mathbb{E}(\zeta^k| \overline{A})$ for $k=1,2$, the first two
moments of
${\zeta}| \overline{A}$. (We will have to do a lot of similar estimations
later.) Note that $0\le\zeta^k \le1$.

This is where we use the last paragraph of Section \ref{sec:sand}.
Recall that $B$ denotes the event that $P(v\ge U^*)\subset\Pi_{\eta}$
and $P(v\le U^*)$ contains at most $3b_2F(P) \ln^d \eta$ points from
$X$. Here, $U=(\ln{\eta})/ \eta$ and $U^*=d6^dU$. Lemma \ref{l:Bpi}
says that $\mathbb{P}(\overline{B})\ll F(P) \eta^{-3d}$. Let $I(B)$
denote the
indicator function of the event $B$. Observe that $\zeta^k I(B) \leq
V (P(v \le U^*))^k$. Moreover,
$V(P(v \le U^*)) \ll F(P)(\ln\eta)^d/ \eta$,
by Theorem \ref{th:G}. So, we have
%
%
\begin{eqnarray}\label{eq:X|A}
\mathbb{E}(\zeta^k| \overline{A})&=& \mathbb{E}\bigl(\zeta^k \bigl(1-I(B)\bigr)|
\overline{A}\bigr)+\mathbb{E}(\zeta^k I(B)| \overline{A})\nonumber\\
&\le&
\mathbb{E}\bigl(\bigl(1-I(B)\bigr)| \overline{A}\bigr)+ V\bigl(P(v \le U^*)\bigr)^k
\\
&\ll& \mathbb{P}(\overline{B} | \overline{A})
+ \biggl(F(P)\frac{\ln^d \eta} \eta\biggr)^k
\ll\biggl(F(P)\frac{\ln^d \eta} \eta\biggr)^k .\nonumber
\end{eqnarray}
Here, we have used the estimate
%
%
\begin{equation}\label{eq:B|A}
\mathbb{P}(\overline{B}|\overline{A})\le
\frac{\mathbb{P}(\overline{B})}{\mathbb{P}(\overline{A}) }
\ll\frac{ F(P) \eta^{-3d}}{(\ln{\eta})^{-(3d)^{d+2}}}
\ll F(P) \eta^{-3d+1},
\end{equation}
where the lower bound for $\mathbb{P}(\overline{A})$ comes from Claim
\ref{cl:PA}.

As for $\mathbb{E}(\zeta^k| A) $, Claim \ref{cl:Aimpl} tells us that
\[
\mathbb{E}(\zeta^k| A) \leq V\bigl(P(v \le T^*)\bigr)^k \ll
\biggl(F(P)\frac{\ln^d \eta} \eta\biggr)^k .
\]
Thus, we get, using Claim \ref{cl:simple}, that
%
%
\begin{equation}\label{eq:kA}
|\mathbb{E}(\zeta^k|A)-\mathbb{E}(\zeta^k)| \ll\biggl(F(P)\frac{\ln^d
\eta} \eta\biggr)^k
\mathbb{P}(\overline{A}).
\end{equation}

We check condition (ii) first. Since $\operatorname{Var}(V(\Pi_{\eta
}))=\operatorname{Var}(1-\zeta)=
\operatorname{Var}(\zeta)=\mathbb{E}(\zeta^2 )-(\mathbb{E}(\zeta
))^2$ and similarly for $\operatorname{Var}
(V(\Pi_{\eta})|A)$, the aim is to estimate
%
%
\begin{eqnarray}\label{eq:var1}
&&\bigl| \bigl(\mathbb{E}(\zeta^2|A) -(\mathbb{E}(\zeta|A))^2\bigr) -
\bigl(\mathbb{E}(\zeta^2)-(\mathbb{E}
(\zeta
))^2\bigr)\bigr| \nonumber\\[-8pt]\\[-8pt]
&&\qquad \le|\mathbb{E}(\zeta^2|A)-\mathbb{E}(\zeta^2)|+|(\mathbb
{E}(\zeta|A))^2-(\mathbb{E}(\zeta))^2|.\nonumber
\end{eqnarray}
The first term in the last line is bounded in (\ref{eq:kA}) with
$k=2$. For the second, we have
%
%
\begin{eqnarray}\label{eq:var2}
|(\mathbb{E}(\zeta|A))^2-(\mathbb{E}(\zeta))^2|
&=&
|\mathbb{E}(\zeta|A)+\mathbb{E}(\zeta)|\cdot|\mathbb{E}(\zeta
|A)-\mathbb{E}(\zeta)|\nonumber\\
&\ll& F(P)\frac{\ln^d \eta} \eta|\mathbb{E}(\zeta|A)-\mathbb
{E}(\zeta)|
\\
&\ll&\biggl(F(P)\frac{\ln^d \eta}
\eta\biggr)^2\mathbb{P}(\overline{A}),\nonumber
\end{eqnarray}
where (\ref{eq:kA}) and (\ref{eq:X|A}) have been applied with $k=1$. We
now need the lower
bound $\operatorname{Var}(V(\Pi_{\eta}))\gg F(P) \eta^{-2}\ln
^{d-1} \eta$ from
Theorem \ref{th:VarPi}.
Combining this lower bound, formulae (\ref{eq:kA}), (\ref{eq:var1}),
(\ref{eq:var2}) and
Claim \ref{cl:PA} yields
%
%
\begin{eqnarray}\label{eq:varVA}
|{\operatorname{Var}(V(\Pi_{\eta})|A)-\operatorname{Var}
(V(\Pi_{\eta}))}|
&\ll&
\biggl(F(P)\frac{\ln^d \eta} \eta\biggr)^2 \mathbb
{P}(\overline{A})
\nonumber\\[-8pt]\\[-8pt]
&\ll& F(P)^2 \frac{\ln^{d+1} \eta}{\ln^{4d^2} \eta}\operatorname
{Var}(V(\Pi
_{\eta})).\nonumber
\end{eqnarray}
This shows that condition (ii) of Lemma \ref{l:transV} is satisfied
with ${\varepsilon}_2({\eta} ) \ll
F(P)^2 \times\break\ln^{-4d^2 +d+1} \eta$. This, together with Theorem \ref
{th:VarPi}, also immediately
proves Claim~\ref{cl:var*}, that is, $\operatorname{Var}(V(\Pi_{\eta
})|A)\gg F(P)
\eta
^{-2} \ln^{d-1} \eta$ when
$\ln^{4d^2 -d-1} \eta\gg F(P)^2$.

Finally, (\ref{eq:kA}) with $k=1$ gives
\begin{eqnarray*}
&&
|\mathbb{E}(V(\Pi_{\eta})|A)- \mathbb{E}(V(\Pi_{\eta
}))| \\
&&\qquad= |\mathbb{E}(\zeta|A)-\mathbb{E}
(\zeta)|
\ll
F(P) \frac{\ln^d \eta} \eta\mathbb{P}(\overline{A})\\
&&\qquad\ll F(P)^2 \frac{\ln^d \eta}{ \eta\ln^{4d^2} \eta}
\ll F(P)^ {3/2} \frac{\ln^{({d+1})/2} \eta}{\ln^{4d^2} \eta}
\sqrt{\operatorname{Var}(V(\Pi_{\eta}))}.
\end{eqnarray*}
Thus, condition (i) is also satisfied with ${\varepsilon}_1({\eta} )
\ll
F(P)^{{3/2}} \ln^{-4d^2 +(d+1)/2} \eta$.

Condition (iii) is the simplest to check: set
$\zeta=I(V(\Pi_{\eta}))\le x)$ and apply Claim~\ref{cl:simple}. Then
\begin{eqnarray*}
|\mathbb{E}(\zeta|A)-\mathbb{E}(\zeta)|&=&\bigl|\mathbb{P}\bigl(V(\Pi_{\eta
})\le x|A\bigr)-\mathbb{P}\bigl(V(\Pi
_{\eta
})\le x\bigr)\bigr| \\
&\le& 2\mathbb{P}(\overline{A})\ll F(P) \ln^{-4d^2} \eta
\end{eqnarray*}
and thus (iii) holds with ${\varepsilon}_3({\eta} ) \ll F(P) \ln
^{-4d^2} \eta$.
\end{pf*}
\begin{pf*}{Proof of Lemma \ref{l:transfs}}
This proof is similar to the
previous one and so we only point out the main differences.
Set $\zeta=f_{\ell}(\Pi_{\eta})$. We want to
estimate, for $k=1,2$,
%
%
\begin{equation}\label{eq:transfs}
\mathbb{E}(\zeta^k|\overline{A})=\mathbb{E}\bigl(\zeta
^k\bigl(1-I(B)\bigr)|\overline{A}\bigr)+
\mathbb{E}(\zeta^kI(B)|\overline{A}).
\end{equation}
Note that, given $B$, $\Pi_{\eta}$ can have at most $ F(P) \ln^d
\eta$ vertices,
implying that $\zeta^k I(B) \ll(F(P) \ln^d \eta)^{k(\ell+1)}$,
which is
an upper bound for the second term in (\ref{eq:transfs}). The first
term needs extra care since the random variable $\zeta$ is not
bounded. Let $N$ be a random variable which is Poisson distributed
with mean $ \eta$ and write $E_m$ for the event $N=m$. Of course,
$\zeta\le m^{\ell+1}\le m^d$ under condition $E_m$. Thus,
\begin{eqnarray*}
&&\mathbb{E}\bigl(\zeta^k\bigl(1-I(B)\bigr)|\overline{A}\bigr) \\
&&\qquad= \sum_{m=0}^{\infty}
\mathbb{E}\bigl(\zeta
^k\bigl(1-I(B)\bigr)|\overline{A} E_m\bigr)\mathbb{P}(E_m)\\
&&\qquad\le \sum_{0 \leq m < 3 \eta} \mathbb{E}\bigl(\zeta
^k\bigl(1-I(B)\bigr)|\overline{A}
E_m\bigr)\mathbb{P}(E_m)\\
&&\qquad\quad{} +\sum_{3 \eta\leq m} \mathbb{E}\bigl(\zeta^k\bigl(1-I(B)\bigr)|\overline{A}
E_m\bigr)\mathbb{P}
(E_m)\\
&&\qquad\le \sum_{0 \leq m < 3 \eta} (3 \eta)^{kd}\mathbb
{E}\bigl(\bigl(1-I(B)\bigr)|\overline
{A} E_m\bigr)\mathbb{P}(E_m)\\
&&\qquad\quad{} +\sum_{3 \eta\leq m } m^{kd}\mathbb{E}\bigl(\bigl(1-I(B)\bigr)|\overline{A}
E_m\bigr)\mathbb{P}
(E_m)\\
&&\qquad\ll (3 \eta)^{kd}\sum_{0 \leq m < 3 \eta}\mathbb{P}(\overline
{B}|\overline
{A} E_m)\mathbb{P}(E_m)+\sum_{3 \eta\leq m }
m^{kd}\mathbb{P}(E_m)\\
&&\qquad\ll (3 \eta)^{kd}\mathbb{P}(\overline{B}|\overline{A} )+\sum_{3
\eta
\leq m }
m^{kd}\mathbb{P}(E_m)\ll F(P) \eta^{-d+1},
\end{eqnarray*}
where we have used (\ref{eq:B|A}) and the routine estimation of
$\sum_{3 \eta}^{\infty} m^{kd} \mathbb{P}(E_m)$ is omitted. Using
this and
Claim \ref{cl:Aimpl} for $\mathbb{E}(\zeta^k|A)$ yields, for $k=1,2$,
\[
\mathbb{E}(\zeta^k|A),\ \mathbb{E}(\zeta^k|\overline{A}) \ll( F(P)
\ln^d \eta
)^{k(\ell
+1)}\leq F(P)^{kd} \ln^{kd^2} \eta.
\]
We again need Claim \ref{cl:simple} and the lower bound from
Theorem \ref{th:VarPi} to show that
%
%
\begin{equation} \label{eq:varfA}
|{\operatorname{Var}}(\zeta|A)- \operatorname{Var}(\zeta)| \ll
F(P)^{2d} \ln^{-2d^2-d+1 } \eta\operatorname{Var}
(\zeta
) .
\end{equation}
So, condition (ii) is satisfied. It also follows that $\operatorname
{Var}(\zeta|A)
\gg
F(P) \ln^{d-1} \eta$ for
$\ln\eta\gg F(P)^{1/d }$,
which is Claim
\ref{cl:varface} from the previous section.

Checking condition (i) follows along the same lines and condition (iii) is
straightforward.
\end{pf*}

\section{Proofs of the auxiliary lemmas}\label{sec:aux}

In this section, we assume that $P$ is a fixed polytope in $\mathbb{R}^d$
whose volume is $1$. We first prove the following claim, where
$\beta=2e d^3 + 1$ (a $\beta$ which is different from the one in
Section \ref{sec:sand}).
\begin{claim} \label{cl:Szs}
For all $T $ and all $z \in P$ satisfying $0< v(z) < \frac12, v(z)
\leq T $, we have
\[
S(z,T)\subset C^{ {\beta T}/{v(z)}} (z).
\]
\end{claim}
\begin{pf}
Set $s=v(z)$. Let $C(z)= \{ x\in\mathbb{R}^d \vert u
\cdot x \geq h_P(u) -h_z \}$ be the minimal cap of $z \in P$ and
denote by $H_h$ the hyperplane $\{ x\in\mathbb{R}^d \vert u \cdot x
= h_P(u) - h \}$. Then $H_0$ touches the boundary of $P$ in a
center of the cap, $H_{h_z}$ is the bounding hyperplane of $C(z)$
and, thus, $z \in H_{h_z}$. Write $Q_{z}= P \cap H_{h_z}$. The
following simple geometric arguments show that for every $ h \in
[0,h_z]$, we have
%
%
\begin{equation} \label{eq:Vd-1}
V_{d-1} (P \cap H_h) \leq2d V_{d-1} (Q_z)
\end{equation}
if $s \leq{\frac12}$, where $V_{d-1}( \cdot)$ stands for
$(d-1)$-dimensional volume. Indeed, the Brunn--Minkowski inequality
shows that for some $h_{\max}$, the volume of the sections $P \cap
H_h$ is first increasing for $h \in[0, h_{\max}]$ and then
decreasing for $h \in[h_{\max}, w]$. Here, $w$ denotes the width of
$P$ in the direction $u$. Thus, if $h_z \in[0, h_{\max}]$, equation
(\ref{eq:Vd-1}) is immediate (with $2d$ replaced by $1$). And if
$h_z > h_{\max}$, then we have to show that $V_{d-1} (P \cap
H_{h_{\max}})
\leq2d V_{d-1} (Q_z)$. Clearly,
\[
\frac1 d w V_{d-1} (P \cap H_{h_{\max}} ) \leq1 .
\]
Since we assume that $h_z > h_{\max}$ here, $V_{d-1} (P \cap H_h)$ is
decreasing for $h \in[h_z, w]$
and we also have
\[
(w-h_z) V_{d-1} (Q_z) \geq1-s \geq\tfrac12 .
\]
Combining this gives
\[
V_{d-1} (P \cap H_{h_{\max}}) \leq\frac d {w} \leq\frac d {w-h_z}
\leq2 d V_{d-1} (Q_z),
\]
which is (\ref{eq:Vd-1}). It follows that
\[
s \leq2d h_z V_{d-1} (Q_z) .
\]

Clearly, the set $S(z,T)$ is the union of caps $C \subset P(v\le T)$
such that $z \in C$. Let $C$ be such a cap. Then $V(C) \le dT$, by
Lemma \ref{l:E}. If $C$ contains a point of $H_h$, then
\[
V(C) \geq\frac1 d (h - h_z) V_{d-1} (C \cap Q_{z}) .
\]
As is well known (see, e.g., \cite{ELR}), $z \in C$ is the
center of gravity of $Q_{z}$. A result of
Gr\"unbaum \cite{grun} then tells us that $V_{d-1} (C \cap Q_{z}) \geq
\frac1 e V_{d-1} (Q_{z})$. Thus,
\[
\frac1 {2 e d^2} \frac{s (h-h_z)}{h_z} \leq\frac1 {e d}
(h-h_z) V_{d-1} (Q_{z}) \leq V(C) \leq dT .
\]
Hence, the distance between an arbitrary point of $S(z,T)$ and $H_0$
is at most $h=\frac Ts 2 e d^3 h_z + h_z \leq(2e d^3 + 1)
\frac{T}{s} h_z$, which shows that, indeed,
$S(z,T)\subset C^{\beta T/s} (z)$.
\end{pf}
\begin{pf*}{Proof of Lemma \ref{l:volS}}
Again, setting $v(z)=s$, the condition is
$0<s \leq\frac12, 2s\leq T$. Choose ${\beta}$ as in
Claim \ref{cl:Szs}. Let $C(z)$ be the minimal cap of $z$ and set
$C^*=C^{\beta T/s}(z)$ and $V^*=V(C^*)$, noting that $C^*$ is a
polytope. By trivial volume estimates, $V^ * \le(\beta T/s)^d
V(C(z))= (\beta T)^d/s^{d-1}$. First, assume that $C^ * =P$ and,
thus, $1/T \leq\beta^ d (T/s)^{d-1}$. Then, since $S(z,T)\subset
P(v \leq T)$, we have
\[
V(S(z,T)) \ll F(P) T \ln^ {d-1} \biggl(\frac{1}{T} \biggr)
\]
for $T \leq s_0$, by Theorem \ref{th:G}, which gives
\[
V(S(z,T)) \ll F(P) T \ln^ {d-1} \biggl(\frac{T} s \biggr)
\]
for any $T$ with $2s \leq T$.
For $C^ * \neq P$, trivial volume estimates show that $V(C^ *) \geq
(\beta T/ds) V(C(z)) =({\beta}/d) T$. Claim \ref{cl:Szs} shows that
\[
S(z,T)\subset P(v_P \leq T) \cap C^* \subset C^*(v_{C^*}\le T),
\]
where we have written $v_{C^*}$ to emphasize that the underlying convex set
is now $C^*$. By Theorem \ref{th:G}, there is a constant $s_0$ such
that for $T \leq s_0 V^ *$,
\begin{eqnarray*}
V(S(z,T)) &\leq&
V\bigl(C^*(v_{C^*}\le T )\bigr) \ll
F(C^*) T \ln^{d-1} \biggl( \frac{V^*}T \biggr) \\
&\le& F(P) T \ln^{d-1} \biggl( \frac{V^*}T \biggr).
\end{eqnarray*}
Here, we used the fact that $F(C^*)\le F(P)$, which can be proven
quite easily (we omit the proof). In the remaining case, $s_0 V^ *
\leq T \leq(d/ {\beta}) V^ *$, we have
\[
V\bigl(C^*(v_{C^*}\le T )\bigr) \leq V^ * \leq s_0^ {-1} T
\ll F(P) T \ln^{d-1} \biggl( \frac{V^*}T \biggr) .
\]
The lemma follows since $V^*/T\le{\beta}^d (T/s)^{d-1}$.
\end{pf*}
\begin{pf*}{Proof of Lemma \ref{l:zinC}}
Assume that $z_i \in Z(s)\cap C$. Then,
by Lemma \ref{l:B}, $M(z_i,1)\subset C^2$.
Thus, for $s \leq s_0$, the set $K'_i(s) = M(z_i,\frac12) \cap C(z_i)$ lies
in $P(v\le s) \cap C^2$. The sets $K'_i(s)$,
$i=1,\ldots,m(s)$, are pairwise disjoint, so the usual volume argument
applies:
\[
|Z(s) \cap C| \ll\frac{V(P(v\le s)\cap C^2)}s
\]
as $V(K'_i(s))\gg s$. Further, $P(v\le s)\cap C^2 \subset
C^2(v_{C^2}\le s)$, whose volume can be estimated in the same way as in
the previous proof. Theorem \ref{th:G} gives
\[
V\bigl(C^2(v_{C^2}\le s)\bigr)\ll F(C^2) s \ln^{d-1} \biggl(\frac{V(C^2)}{s} \biggr) \ll F(P)s
\ln^{d-1} \biggl(\frac Ts \biggr)
\]
for $s \leq s_0 V(C^2)$ since $F(C^2)\ll F(P)$.
And, for $s_0 V(C^2) \leq s \leq s_0 $, the lemma follows from the fact that
$V(C^2(v_{C^2}\le s)) \leq V(C^2)$ and $s \leq2T$.
\end{pf*}
\begin{pf*}{Proof of Lemma \ref{l:zinK}}
Since $V(K_j^{\lambda}(T)) \le {\lambda}^d
V(K_j(T))\le(6{\lambda})^dT$, each $y_j \in Y(T)$ with $z \in
K_j^{\lambda}(T)$ is
contained in
$S(z,(6{\lambda})^dT)$. It is also clear that $M(y_j,{\frac12}) \cap C(y_j)$
lies in
$S(z,(6{\lambda})^dT)$, once $y_j \in Y(T)$. Thus, the usual volume argument
applies, with the upper bound on $V(S(z,(6{\lambda})^dT))$ coming from
Lem\-ma~\ref{l:volS}.
\end{pf*}
\begin{pf*}{Proof of Lemma \ref{l:nkl}}
Let $C(z)$ be the minimal cap of $z$.
Claim \ref{cl:Szs} shows that $S(z,T)$ is contained in the
cap $C:=C^{\beta T/s} (z)$ with volume $ V(C) \leq({\beta}
T)^d/s^{d-1}$. Lemma \ref{l:zinC} then applies and gives
\[
|Z(s)\cap S(z,T)|\le|Z(s)\cap C| \ll F(P) \ln^{d-1}\frac{V(C)} s
\ll F(P) \ln^{d-1} \biggl(\frac Ts \biggr)
\]
for $s \leq s_0$ and $2s \leq V(C)$ since $V(C)/s \ll(T/s)^d$. The
inequality $2s \leq V(C)$ follows from the trivial volume estimate
if $C \neq P$, and from $s \leq s_0$ if $C=P$.
\end{pf*}
\begin{pf*}{Proof of Claim \ref{cl:convxi1}}
Clearly, it suffices to show that each
cap $C$ whose bounding hyperplane touches $P(v
\geq T^*)$ contains at least one point $x_i$. If this is not the
case, then there is a cap $C$ whose bounding hyperplane touches
$P(v\ge T^*)$ with no $x_i \in C$ and thus no $M(y_i, {\frac12})
\subset
C$ either.

We now claim that $C^{1/3}$ is disjoint from all Macbeath
regions $M(y_i, {\frac12})$. Assume, for simpler notation, that $u
\cdot
x=h$ with $h>0$ is the equation of the bounding hyperplane of $C$,
and $u \cdot x=0$ is the equation of the supporting hyperplane of $P$
and $C$. If $u
\cdot y_i=g$, then $M(y_i,1)$ lies between hyperplanes $u \cdot
x=2g$ and $u \cdot x=0$. Thus, $M(y_i,{\frac12})$ lies between hyperplanes
$u \cdot x=\frac32 g$ and $u \cdot x=\frac12 g$. Here, $\frac32 g
>h$ holds since, otherwise,
$M(y_i,{\frac12}) \subset C$. Then $g > \frac23 h$, implying that $u
\cdot
x=\frac13 h$ is a separating hyperplane between $M(y_i,{\frac12})$ and
$C^{1/3}$. This proves the claim.

By trivial volume estimates, $V(C^{1/6})$ is at least $dT$. Let
$x_0$ be the point in $C^{1/6}$ where $v(x)$ takes its maximal
value on $C^{1/6}$. By Lemma \ref{l:E}, $V(C^{1/6}) \le
dv(x_0)$ and so $v(x_0) \ge T$. This shows the existence of a point
$z \in P(v=T) \cap C^{1/6}$. However, we then have that $M(z,
{\frac12}) \subset
C^{1/3}$ is disjoint from all $M(y_i, {\frac12})$, which is impossible
since $Y(T)=\{y_1, \ldots, y_{m(T)}\}$ is a saturated system.
\end{pf*}

\section*{Acknowledgments}
The authors are indebted to an anonymous referee who pointed out
several inconsistencies in an earlier version of this paper. The
present, thoroughly rewritten and hopefully more readable version owes
much to this referee. Part of this paper was written during a pleasant
and fruitful visit, by the first author, to the Institute for Advanced
Study at the Hebrew University of Jerusalem, where the excellent
ambiance and working conditions were greatly appreciated.

%

%
\printaddresses

\end{document}